\documentclass[amscd,amssymb,verbatim,12pt]{amsart}
\usepackage{epsfig}
\usepackage{graphicx}
\usepackage{pdfpages}

\newif\ifAMS
\IfFileExists{amssymb.sty}
  {\AMStrue\usepackage{amssymb}}
  {\usepackage{latexsym}}
  \usepackage{enumerate}
\theoremstyle{plain}

\newtheorem*{First}{Theorem \ref{T:Main}}
\newtheorem*{Second}{Theorem \ref{second}}
\newtheorem*{Ricci}{Corollary \ref{ricci}}

\newtheorem{Thm}{Theorem}[section]
\newtheorem{Cor}[Thm]{Corollary}
\newtheorem{Lem}[Thm]{Lemma}

\theoremstyle{definition}
\newtheorem{Def}{Definition}

\theoremstyle{remark}

\newtheorem{Qu}[Thm]{Question}

\newtheorem{Ex}[Thm]{Example}

\DeclareMathOperator{\asdim}{asdim}
\newcommand{\interior}{^{ \kern-5pt ^\circ}}

\begin{document}
\title
{Polynomial growth and asymptotic dimension}

\author
{Panos Papasoglu}

%\subjclass{53C23}

\email {} \email {papazoglou@maths.ox.ac.uk}

\address
{Mathematical Institute
University of Oxford
Andrew Wiles Building
Radcliffe Observatory Quarter (550)
Woodstock Road
Oxford
OX2 6GG }

\begin{abstract} 
Bonamy et al \cite{BBEGLPS} showed that graphs of polynomial growth have finite asymptotic dimension.
We refine their result showing that a graph of polynomial growth strictly less than $n^{k+1}$ has
asymptotic dimension at most $k$. As a corollary Riemannian manifolds of bounded geometry and polynomial 
growth strictly less than $n^{k+1}$ have asymptotic dimension at most $k$.

We show also that there are graphs of growth $<n^{1+\epsilon}$ for any $\epsilon >0$ and infinite asymptotic Assouad-Nagata dimension.

\end{abstract}
\maketitle

\section{Introduction}

Asymptotic dimension is a large scale analog of topological dimension that was introduced by Gromov \cite{Gr}.
It is invariant under quasi-isometries and even stronger under coarse embeddings, so one can think of it as a
large scale topological notion (see \cite{BD} for an introduction to the subject).

Asymptotic dimension is relevant in several contexts: in geometric group
theory, as groups of finite asymptotic dimension satisfy the Novikov conjecture \cite{Y}, 
in geometry \cite{LS}, \cite {BL}
%where it is conjectured \cite{Gr2} that universal covers of $n$-manifolds of strictly positive scalar curvature have asymptoticdimension bounded by $n-2$
and in graph theory \cite{OR}, \cite{FP},\cite{BBEGLPS}.

The \textit{asymptotic dimension} ${\rm asdim}\,X$ of a metric space $X$ is defined as follows: ${\rm asdim}\,X \leq n$ if and only if for every $m > 0$ there exists $D(m)>0$ and a covering $\mathcal{U}$ of $X$ by sets of diameter $\leq D(m)$ 
($D(m)$-bounded sets) such that any $m$-ball in $X$ intersects at most $n+1$ elements of $\mathcal{U}$.

We say that the $m$-multiplicity of the cover $\mathcal{U}$ is at most $n+1$.

If we can pick $D(m)$ to be a linear function in the definition above, then we say that the \textit{asymptotic Assouad-Nagata dimension} of $X$ is bounded by $n$ (\cite{DS},\cite{BNLM}).

\v{S}pakula and Tikuisis, ask in \cite{ST} (see the footnote in page 3) whether spaces of polynomial growth have finite asymptotic dimension
(and apparently as this is a quite natural question, it was considered by other people as well).

A different notion of dimension was considered earlier by Linial, London, Rabinovich \cite{LLR}, and Linial \cite{Lin}, namely one defines the
dimension of a graph $G$ to be the smallest $n$ for which there is an embedding $f:G\to \mathbb R^n$ so that $d(f(u),f(v))\geq 1$ for $u\ne v$ and for some $c>0$,
$d(f(u),f(v))\leq c$ if $u,v$ are adjacent. 
Krauthgamer and Lee \cite{KL} showed that graphs of polynomial growth $\gamma (r)\leq Cr^k$ embed in this sense in $\mathbb R^{O(k\log k)}$.

Using the result of \cite{KL} Bonamy et al prove in \cite{BBEGLPS} that graphs of polynomial growth $\leq Cr^k$ have asymptotic dimension
bounded by $O(k\log k)$ answering the question of \v{S}pakula-Tikuisis. It is further shown that graphs of superpolynomial growth can have infinite asymptotic dimension.

Benjamini and Georgakopoulos \cite{BG} show that planar triangulations of subquadratic growth are quasi-isometric to trees (and so they
have $\asdim =1$).

From the geometric point of view it is interesting to calculate the exact asymptotic dimension of a space. This has been accomplished for several `natural' classes of spaces: It is shown in \cite{BL} that the asymptotic dimension of a hyperbolic group $G$ is equal to $\dim (\partial G)+1$, and in \cite{Gr} that $n$-dimensional Hadamard manifolds of pinched negative curvature have asymptotic dimension $n$ (see also \cite{LS} for a detailed proof and and an extension of this to asymptotic Assouad-Nagata dimension). It is shown in \cite{FP}, \cite{JL},\cite{BBEGLPS} that planar graphs (or more generally planar geodesic metric spaces) have asymptotic dimension at most 2. In this paper we extend this list to the class of spaces with polynomial growth.
Note that spaces with polynomial growth appear in several settings. For example doubling spaces have polynomial growth \cite{He} and manifolds of non-negative Ricci curvature have polynomial volume growth \cite{Gr2}. It is shown in \cite{LS} that doubling metric spaces have finite asymptotic dimension (in fact also finite Nagata dimension)
and a sharp bound of their Nagata dimension (hence also asymptotic dimension) in terms of Assouad dimension is given in \cite{DR}.
Tessera in \cite{T} studies geometric properties of general graphs of polynomial growth.

We state now our results.
We view a connected graph as a geodesic metric space
where each edge has length $1$. 

\begin{Def}\label{growth}
We define the \textit{growth function} of a connected graph $G=(V,E)$ to be
$$\gamma (r)=\sup \{|B_v(r)|:v\in V\}.$$
\end{Def}
Where we denote above by $|X|$ the number of vertices of a subset $X$ of $G$. Due to the independence from a base vertex, some authors call $\gamma (r)$ uniform growth function.

We prove the following:

\begin{First}
Let $G=(V,E)$ be a connected graph with growth function $\gamma (r)$ satisfying 
$$\lim _{r\to \infty }\dfrac {\gamma (r)}{r^{k+1}}=0\ \ \ \text { for some }k\in \mathbb N.$$
Then ${\rm asdim}\,G\leq k$.

\end{First}

%\kf{is it uniformly at most two? ie, uniform for all metrics on a plane?}

As a corollary we have:

\begin{Ricci}
If $M^n$ is a Riemannian manifold of bounded geometry and volume growth function $Vol(r)$ satisfying 
$$\lim _{r\to \infty }\dfrac {Vol(r)}{r^{k+1}}=0\ \ \ \text { for some }k\in \mathbb N $$
then ${\rm asdim}\, M^n\leq k$.

\end{Ricci}

We define the volume growth function for $M^n$ as for graphs
$$Vol(r)=\sup \{vol(B_x(r)):x\in M^n\}.$$
We say that $M^n$ is of bounded geometry if there are $a>0,b>1$ such that for any open ball of radius $a$ in $M^n$
there is a $b$-bilipschitz map to a Euclidean open ball of radius 1. Recall that a map $f:X\to Y$ is $b$-bilipschitz if it is onto and
$1/b\, d(x,y)\leq d(f(x),f(y))\leq b\, d(x,y)$ for all $x,y\in X$. We remark that other common definitions of bounded geometry
for non-compact manifolds imply ours.

It turns out that Theorem \ref{T:Main} applies more generally to metric spaces for an appropriate definition of volume (see section 3).
However one does not have a similar bound for the asymptotic Assouad Nagata dimension. We have:

\begin{Second}
There is a metric space $X$ with 1-growth function satisfying 
$$\lim _{r\to \infty}\dfrac {\gamma (r)}{r}=0$$ and infinite asymptotic Assouad-Nagata dimension.

There is a graph $(G,E)$ with growth function $\gamma (r)$ 
satisfying 
$$\lim _{r\to \infty}\dfrac {\gamma (r)}{r^{1+\epsilon}}=0$$ 
for any $\epsilon >0$ and infinite asymptotic Assouad-Nagata dimension.

\end{Second}

We note it follows from this theorem that the bound of asymptotic dimension in terms of Assouad dimension implied by \cite{DR} is far from optimal. Since the Assouad dimension bounds the Nagata dimension, for the graphs of the theorem the Assouad dimension
is infinite while the asymptotic dimension is equal to 1.

 \section*{Acknowledgements}
I thank  Agelos Georgakopoulos for interesting discussions, I thank Urs Lang, Romain Tessera and Panagiotis Tselekidis for their comments on a first draft of this paper and Alexander Engel for bringing \cite{ST} to my attention.
 %\kf{add more}

\section{Graphs of polynomial growth.}

\begin{Def}
Let $G=(V,E)$ be a connected graph and $r>0$. We say that $X\subset G$ is \textit{$r$-scale connected} if
for any $x,y\in X$ there is a sequence  $x_1=x,x_2,...,x_n=y$ in $G$ such that $d(x_i,x_{i+1})\leq r$
for all $r$. If $A\subset X$ we say that $A$ is an \textit{$r$-connected component} of $X$ if $A$ is a maximal
$r$-scale connected subset of $X$.
\end{Def}

So for example $V$ has a single 1-connected component (equal to itself) and each singleton $\{v\}$ is a $1/2$-connected component of $V$.

\begin{Def}
Let $G=(V,E)$ be a connected graph and $r>0$. We say that $X\subset G$ has \textit{$r-\dim X\leq n$} if
there is a $D>0$ and a cover $\mathcal U$ of $X$ by sets of diameter $\leq D$ such that any $r$-ball intersects at 
most $n+1$ elements of $\mathcal U$. 
\end{Def}

Clearly $\asdim G\leq n$ if $r-\dim G\leq n$ for all $r$.

\begin{Def}
Let $G=(V,E)$ be a connected graph and let $X\subset G$. We define the \textit{growth function of $X$} to be
$$\gamma_X (r)=\sup \{|B_v(r)\cap X|:v\in V\}.$$
\end{Def}

\begin{Lem}\label{zero}
Let $G=(V,E)$ be a connected graph and let $X\subset G$. If there is some $n$ so that $\gamma_X (n)<\dfrac{1}{2r} n$ 
then $r-\dim X=0$.
\end{Lem}

\proof
It suffices to show that the $r$-connected components of $X$ are bounded. By hypothesis there
is an $n$ such that $\gamma_X (n)<n/2r$. It follows that any $r$-connected component of $X$ has diameter bounded by $n$.
\qed

\begin{Def}
Let $G=(V,E)$ be a connected graph and let $X\subset V$. We say that $X$ is \textit{$(D,r)$-separating} if
all $r$-connected components of $V\setminus X$ have diameter bounded by $D$.

More generally if $Y\subset V$ and $X\subset Y$ we say that that $X$ is a \textit{$(D,r)$-separating} subset of $Y$ if
all $r$-connected components of $Y\setminus X$ have diameter bounded by $D$.
\end{Def}

We fix now a vertex $e\in G$. 

\begin{Def}
We say that a $(D,r)$-separating subset $X$ of $V$ is \textit{minimal} if there is a sequence of $(D,r)$-separating subsets $X_n$ of $V$
such that the following hold:

\item for any $(D,r)$-separating
subset $Z$ of $Y$
$$|B_e(n)\cap Z|\geq |B_e(n)\cap X_n| $$

\item for any $k>0$ there is some $n_k\geq k$ such that
$$B_e(k)\cap X=B_e(k)\cap X_{n_{k+t}} \text{ for all } t\in \mathbb N.$$

%for any $k>0$ there is a 
 %for any other $(D,r)$-separating
%subset $Z$
%$$|B_e(n)\cap Z|\geq |B_e(n)\cap X|\,\, \text{ for all } n\in \mathbb N.$$
%
More generally if $X\subset Y\subset V$ and $X$ is a $(D,r)$-separating subset of $Y$ we say that $X$ is minimal
if the same two conditions are satisfied for $(D,r)$-separating subsets of $Y$.

%for any other $(D,r)$-separating
%subset $Z$ of $Y$
%$$|B_e(n)\cap Z|\geq |B_e(n)\cap X| \,\, \text{ for all } n\in \mathbb N.$$
\end{Def}

\begin{Lem} Let $G=(V,E)$ be a locally finite connected graph.
For any $D,r>0$ minimal $(D,r)$-separating subsets exist. The same is true
for minimal $(D,r)$-separating subsets of a subset $Y$.
\end{Lem}

\proof
Clearly $(D,r)$-separating sets exist, e.g. take the complement of a $D/2$-ball.
Let $X_n$ be a $(D,r)$-separating subset of $V$ for which
$|B_e(n)\cap X_n|$ attains the minimal value among all $(D,r)$-separating subsets.
Since $G$ is locally finite we can pass to a subsequence $X_{n_k}$ such that
$$X_{n_k}\cap B_e(k)=X_{n_{k+t}}\cap B_e(k)$$
for all $t\in \mathbb N$. Then set 
$$X=\bigcup _{k=1}^\infty (X_{n_k}\cap B_e(k)).$$
Clearly $X$ is a minimal $(D,r)$-separating subset of $V$.

The same proof applies for  $(D,r)$-separating subset of subsets of $V$.

\qed

\begin{Def}
Let $G=(V,E)$ be a graph.
An \textit{annulus} with center $v\in G$ and radii $m<n$ is the set
$$A(v,m,n)=\{x\in V: m\leq d(v,x)\leq n\}.$$
We say that $n-m$ is the \textit{thickness} of the annulus.
\end{Def}

\begin{Lem}\label{induction}
Let $k,r\in \mathbb N$ and
let $G=(V,E)$ be a connected locally finite graph. If $X\subset G$ is such that for some $m$ the growth at $n_0=4^mr$ satisfies
$$\gamma _X(n_0)\leq \left(\dfrac {1}{4r}\right) ^{k}  n_0^{k},$$ 
then $r-\dim (X)\leq k-1$.
\end{Lem}

\proof By lemma \ref{zero} the assertion holds for $k=1$. Assume inductively that the assertion holds for $k-1$,
and let
$X$ be as in the lemma. Let $D=2n_0$ and let $Y$ be a minimal
$(D,r)$-separating subset of $X$. If $n_1=n_0/4$ we claim that 
$$\gamma _Y(n_1)\leq \left(\dfrac {1}{4r}\right)^{k-1}  n_1^{k-1}.$$

Indeed assume that this is not the case, so there is some $v$ such that $$|B_v(n_1)\cap Y|>\left(\dfrac{1}{4r}\right)^{k-1}  n_1^{k-1}.$$
Consider the annulus $A(v,n_1,n_0)$. Since $n_0=4n_1$, $A(v,n_1,n_0)$ contains $n_0/2r$ disjoint annuli $A_i$ of thickness $r$.
Clearly $$\sum |A_i|\leq \left(\dfrac {1}{4r}\right)^{k}  n_0^{k}$$ so some annulus, say $A_j$,
satisfies $$|A_j|\leq  \dfrac {2r}{n_0}\left(\dfrac {1}{4r}\right)^{k}  n_0^{k}<\left(\dfrac {1}{4r}\right)^{k-1}  n_1^{k-1}.$$ 

Let $N=d(e,v)+4n_0$.

Since $Y$ is minimal there is some $(D,r)$-separating subset $Y'$ of $X$ such that:
$$B_v(4n_0)\cap Y=B_v(4n_0)\cap Y'$$
and for any $(D,r)$-separating subset $Z$ of $X$
$$|B_e(N)\cap Z|\geq |B_e(N)\cap Y'|$$

Consider now the set $Y_1=(Y'\setminus B_v(n_1))\cup A_j$. It is easy to see that $Y_1$ is also a $(D,r)$-separating subset of $X$.
This contradicts the second property of $Y'$ since clearly
$$|B_e(N)\cap Y_1|<|B_e(N)\cap Y'|.$$

This proves the claim. Note now that $Y$ satisfies the inductive hypothesis so $r-\dim Y\leq k-2$.
Take now a uniformly bounded cover $\mathcal V$ of $Y$ of $r$-multiplicity $\leq k-1$ and add to it the $r$-connected components
of $X\setminus Y$, that have diameter at most $D$ by hypothesis. We obtain a uniformly bounded cover $\mathcal U$ of $X$ of $r$-multiplicity $\leq k$. We note that the diameter of the sets in the cover is bounded by $2n_0$. In particular the diameter
of the sets in the cover that we constructed depends only on the function $\gamma _X$ (and $r$) and not on $X$.

\qed

\begin{Thm}\label{T:Main}
Let $G=(V,E)$ be a connected graph with growth function $\gamma (r)$ satisfying 
$$\lim _{r\to \infty }\dfrac {\gamma (r)}{r^{k+1}}=0\ \ \ \text { for some }k\in \mathbb N.$$
Then ${\rm asdim}\,G\leq k$.
\end{Thm}

%\begin{First}
%Let $G=(V,E)$ be a graph with growth function $\gamma (r)$ satisfying 
%$$\lim _{r\to \infty }\dfrac {\gamma (r)}{r^{k+1}}=0\ \ \ \text { for some }k\in \mathbb N.$$
%Then ${\rm asdim}\,G\leq k$.
%
%\end{First}
%

\proof
By assumption for any $r>0$ there is an $m$ such that for $n_0=4^mr$ 
$$\gamma (n_0)\leq \left(\dfrac {1}{4r}\right) ^{k+1}  n_0^{k+1}.$$ 
So by lemma \ref{induction} $r-\dim G\leq k$. Since this holds for every $r$, $\asdim G\leq k$.

\qed
\begin{Ex} The Cayley graph of $\mathbb Z ^k$ has growth $\gamma (r)\sim r^k$. The above theorem implies that
there is no graph of growth between $r^k$ and $r^{k+1}$ that has asymptotic dimension $k+1$-ie there is nothing `smaller' that
the Cayley graph of $\mathbb Z^{k+1}$ with $\asdim =k+1$. It is important of course for this that we define growth `uniformly' independently of a base vertex. If one defines growth with respect to a base point, paraboloids of dimension $k$ have strictly smaller growth than $\mathbb R^k$ but asymptotic dimension $=k$.

\end{Ex} 
\begin{Ex} \label{ptree}
It is clear that there is no lower bound on $\asdim $ in terms of volume. For example 3-regular trees have exponential
volume growth and $\asdim =1$. It is easy to create similar examples of trees with polynomial growth as well by taking a sparse
set of branch points.

For example take an infinite complete binary tree $T$, with root $e$ (so every vertex has two `children') and subdivide edges so that an edge at distance $n$ from $e$ in $T$ has length $2^{\lfloor n/k\rfloor}$ after the subdivision, for some fixed $k$. Then 

$$\gamma (r)\sim r^{k+1}.$$

\end{Ex}

\section{Assouad-Nagata dimension and metric spaces.}

There are several ways to assign volume functions to general metric spaces. A quite naive definition
is appropriate for this paper:

\begin{Def}\label{net}
Let $X$ be a metric space and $\epsilon, \delta >0$. $N$ is an $(\epsilon, \delta )$-net of $X$ if for
any $x\in X$ there is $n\in N$ with $d(x,n)\leq \epsilon $ and for any $n_1,n_2\in N\, d(n_1,n_2)\geq \delta$.
\end{Def}

Using Zorn's lemma it is easy to see that any metric space $X$ contains an $(\epsilon, \epsilon )$-net for any $\epsilon >0$.

\begin{Def}\label{growth-metric}
Let $X$ be a metric space and let $t\in \mathbb R$. 
We define the $t$-\textit{growth function} of a metric space to be
$$\gamma ^t (r)=\sup \{|B_v(r)\cap N|: v\in X, N \text{  }(t,t)-\text{ net}\}.$$

We say that $X$ has \textit{polynomial growth} if there are $t>0,k,C>0$ such that $\gamma ^t (r)<Cr^k$ for all $r$.
\end{Def}

The supremum in the definition above is over all $(t,t)$-nets (and $v\in X$). 

Observe that $\gamma ^t (r)\leq \gamma ^s (r)$ if $t>s$. Note also that if $\gamma ^t (r)\ne \infty$ for any $r$ then for any $s>t$ there is some $C_s>0$ such that $\gamma ^t (r)\leq C_s \gamma ^s (r)$. Indeed this is because there is a $C>0$ such that if $N$ is a $(t,t)$-net any ball of radius $2s$ contains at most $C$ points of $N$. 

It follows that if $\gamma ^t (r)\ne \infty$ for any $r$ for some $t$ then the 
specific `scale' $s>t$ that we pick to define the growth function does not affect much the asymptotic behavior of the growth function.
For these reasons we will omit below the reference to scale and we will denote the growth function for a specific scale simply by $\gamma (r)$.

Similarly one sees that the growth does not depend much on the net we pick. If $\gamma ^t (r)\ne \infty$ for all $r$ and $N_1,N_2$ are
$(t,t)$-nets then there is $C>0$ such that $|B_v(r)\cap N_1|\leq C|B_v(r)\cap N_2|$.

Note however that if $\gamma ^t (r)= \infty$ for some $r$, it is possible that $|B_v(r)\cap N|<\infty $ for some specific net $N$:

\begin{Ex} Consider the linear graph with vertices $\mathbb N$ where we attach $n$ extra edges to the vertex $n$ for every $n$. Then 
$\gamma ^1 (1)= \infty$ but if we take $\mathbb N$ as a (1,1)-net $|B_v(1)\cap \mathbb N|\leq 3 $ for all $v$.
\end{Ex}
We note that for a graph $\gamma (r)$ is finite for all $r>0$ if and only if it has uniformly bounded degree if and only if $\gamma (r)$ is finite for some $r>0$.

\begin{Lem} \label{metric} Let $X$ be a metric space with $1$-\textit{growth function} $\gamma (r)$ satisfying
$$\lim _{r\to \infty }\dfrac {\gamma (r)}{r^{k+1}}=0\ \ \ \text { for some }k\in \mathbb N.$$
Then ${\rm asdim}\,X\leq k$.
\end{Lem}

\proof

Let $N$ be a $(1,1)$-net of $X$. Given $m>0$ we define a graph $G=(N,E)$ where $(x,y)\in E$ if and only if $d(x,y)\leq 2m$.
$G$ could have several connected components, however for each connected component $\Gamma$ we have that
its growth function satisfies $$\gamma _{\Gamma}(r)\leq \gamma (2mr)$$ so
$$\lim _{r\to \infty }\dfrac {\gamma _{\Gamma} (r)}{r^{k+1}}=0.$$

By theorem \ref{T:Main} the asymptotic dimension of $\Gamma $ is at most $k$, so the $m-\dim$ of $\Gamma $ is bounded by $k$.
In fact as we showed in lemma \ref{induction} there is a cover $\mathcal U$ of $\Gamma $ with $m$-multiplicity $\leq k+1$ such
that the diameter of every sets in $\mathcal U$ is bounded by a constant that depends only on $\gamma $ (see the last line
of the proof of lemma \ref{induction}).

Since this is true for every connected component of $G$ we have that $m-\dim$ of $X$ is bounded by $k$. 

As this is true for every $m$, $\asdim X\leq k$.

\qed

\begin{Cor}\label{ricci}
If $M^n$ is a Riemannian manifold of bounded geometry and volume growth function $Vol(r)$ satisfying 
$$\lim _{r\to \infty }\dfrac {Vol(r)}{r^{k+1}}=0\ \ \ \text { for some }k\in \mathbb N $$
then ${\rm asdim}\, M^n\leq k$.

\end{Cor}
\proof

Let's say that any $a$-ball $B_x(a)$ of $M^n$ is $b$-bilipschitz with the 1-ball of $\mathbb R^n$.
Then $$c_1\leq vol B_x(a)\leq c_2$$
for some $c_1,c_2>0$ and for any $x\in M^n$.

If $N$ is an $(a,a)$-net of $M^n$ and $B_v(r+a)$ is a ball of radius $r+a$ we consider
the set $C=B_v(r)\cap N$. The open balls $B_c(a)$ where $c\in C$ are disjoint, each has volume 
$\geq c_1$ and they are all contained in $B_v(r+a)$. It follows that $c_1\cdot |C|\leq Vol(B_v(r+a))$.

Hence if we see $M^n$ as a metric space its $a$-growth function satisfies
$$\gamma ^{a}(r)\leq \dfrac {Vol(r+a)}{c_1}$$
so by lemma \ref{metric} $\asdim M^n\leq k$.

\qed

\begin{Ex} The bounded geometry hypothesis is necessary: Consider any graph $G$ of bounded degree with infinite asymptotic dimension.
`Thicken' the graph to a 2-manifold $S$ (so edges become thin cylinders). By replacing the edges by thinner and thinner `tubes'
as we go to infinity we obtain a 2-manifold of finite area and infinite asymptotic dimension.
\end{Ex}

%\begin{Lem} 
%\end{Lem}
%
%\proof
%
%\qed

\begin{Thm}\label{second}
There is a metric space $X$ with 1-growth function $\gamma (r)$ satisfying 
$$\lim _{r\to \infty}\dfrac {\gamma (r)}{r}=0$$ and infinite asymptotic Assouad-Nagata dimension.

There is a graph $G=(V,E)$ with growth function $\gamma (r)$ 
satisfying 
$$\lim _{r\to \infty}\dfrac {\gamma (r)}{r^{1+\epsilon}}=0$$ 
for any $\epsilon >0$ and infinite asymptotic Assouad-Nagata dimension.

\end{Thm}

\proof
We give first a sketch of this proof. Both parts are similar so we explain the idea for
the case of graphs, which is more involved. To ensure that $G$ has infinite asymptotic Assouad-Nagata dimension
it suffices for $G$ to contain for each $n$ subgraphs isomorphic to balls of radius $n$ of the Cayley graph of $\mathbb Z^n$.
Of course such balls would lead to big growth, however this can be corrected by rescaling the length of the edges
(or equivalently subdividing each edge many times so that most vertices have degree 2 after this subdivision).
There remains the problem that the degree of the original vertices (before subdivisions) is $n$ (so unbounded) but this can be mended
by replacing the vertices with trees with sublinear growth and $n$-end points.

We proceed now with the details.
Let $\Gamma _n$ be the Cayley graph of $\mathbb Z^n$ with respect to the standard generating set.
Let $B_e(n)$ the ball of radius $n$ in $\Gamma _n$ with center $e$ (so $B_e(n)$ is an $n$-`cube' of side length $2n$). Let $G_n$ be the metric space obtained by $B_e(n)$ by changing the length of edges from 1 to $2^{n^2}$
(so we rescale the metric by $2^{n^2}$). We consider a metric on $$\bigcup _{n\in \mathbb N} G_n$$
so that the $G_n$'s are far apart. For example we may consider the linear graph with vertex set $\mathbb N$ and identify
a vertex of $G_n$ with $10^{n^2}$. We get a graph $\Gamma $ which we see as a metric space containing $\bigcup _{n\in \mathbb N} G_n$ and so we obtain
an induced metric on $\bigcup _{n\in \mathbb N} G_n$. If $X$ is the vertex set of $\bigcup _{n\in \mathbb N} G_n$ for this metric
then the 1-growth function $\gamma (r)$ satisfies 
$$\lim _{r\to \infty}\dfrac {\gamma (r)}{r}=0.$$
As $X$ contains bigger and bigger copies of the vertex sets of $G_n$'s it is easy to see that the asymptotic Assouad-Nagata dimension of $X$ is infinite.
Indeed suppose that the Assouad-Nagata dimension of $X$ is equal to $k$. Then there is a $C>0$ such that for any sufficiently big $r$ there is a cover of $X$ by sets of diameter
$\leq Cr$ such that any $r$-ball intersects at most $k+1$ of these sets. Let's denote by $V(n)$ the vertex set of the $n\times n$ grid in $\Gamma _{k+1}$ 
and let $d_{n}$ be the induced metric on $V(n)$ by the inclusion in $\Gamma _{k+1}$,
then it is clear that $X$ contains copies of $(V(n),2^{n^2}d_n)$ for any $n>0$. So for any sufficiently big $r$ there is a cover $\mathcal U _n$ of $(V(n),2^{n^2}d_n)$ by sets of diameter
$\leq Cr2^{n^2}$ such that any $2^{n^2}r$-ball intersects at most $k+1$ of these sets. By rescaling the metric we have that $(V(n),d_n)$ has a cover $\mathcal U _n'$
by sets of diameter $\leq Cr$ such that any $r$-ball intersects at most $k+1$ of these sets. However as this is true for any $n$ it implies that the asymptotic dimension 
of $\Gamma _{k+1}$ is at most $k$, a contradiction. This proves the first part of the theorem.

To prove the second part we make a similar construction. Let $B_e(n,k)$ be the ball of radius $k$ in $\Gamma _n$ (the Cayley graph of $\mathbb Z^n$) Let $G(n,k)$ be the metric space obtained by $B_e(n,k)$ by changing the length of edges from 1 to $2^k$
(so we rescale the metric by $2^k$). Clearly we can turn this into a graph by subdividing the original edges into $2^k$ edges. We consider the linear graph with vertex set $\mathbb N$ and identify
a vertex of $G(n,k)$ with $10^k$. We obtain in this way a graph $Y_n$ with asymptotic Assouad Nagata dimension at least $n$ and growth $\gamma (r)$ that satisfies 
 $$\lim _{r\to \infty} \dfrac {\gamma (r)}{r^{1+\epsilon}}=0$$ for any $\epsilon >0$. We note that the maximum degree of a vertex in $Y_n$ is $2n+2$ so it is not
possible to obtain a graph of bounded degree taking a `union' of $Y_n$'s.  For this reason the growth function of such a union
would be infinite for any $r\geq 1$.

To correct this we replace each vertex $v$ of degree $2n>2$ of $B_e(n,k)$ in $Y_n$ by a finite binary tree $T_n$ with diameter $2^n$ and $2n$ end vertices. Then we identify each end-vertex of $T_n$ with a vertex of an edge adjacent to $v$.
We call the graph obtained by replacing vertices in this way by $X_n$.

More specifically $T_n$ is a finite tree with $\sim \log n$ branch points and edges of length $\sim 2^n/\log(\log n)$. As usual
we subdivide the edges of $T_n$ into edges of length 1 to get a simplicial graph.

The growth function $\gamma _n$ of $X_n$ satisfies 
$$\lim _{r\to \infty}\dfrac {\gamma _n (r)}{r^{1+\epsilon}}=0$$ 
for any $\epsilon >0$. We claim that $X_n$ has Assouad Nagata dimension at least equal to $n$.
Indeed we have a map $f:X_n\to Y_n $ obtained by collapsing the $T_n$'s to points. Clearly $f$ is a quasi-isometry as
the fibers of this map have uniformly bounded
diameter, so $X_n$ and $Y_n$ have the same asymptotic Assouad Nagata dimension.

Finally we may take $X$ to be the `union' of $X_n$'s.

More precisely take $L$ to be the linear graph with vertex set $n$ and identify
a vertex of $X_n$ with the vertex $2^{2^n}$ of $L$. Each vertex of the graph $X$ obtained in this way has bounded degree.
We note that if we consider the union of $T_n$'s
its growth is bounded by $10r^{1+\epsilon}$ for any $\epsilon >0$ and any $r\geq 1$. It is easy to see that
the growth function $\gamma (r)$ of $X$ satisfies 
$$\lim _{r\to \infty}\dfrac {\gamma (r)}{r^{1+\epsilon}}=0$$ 
and $X$ has infinite asymptotic Assouad Nagata dimension since it contains a copy of $X_n$ for all $n$.

\qed

\section{Discussion and questions}

We say that a proper metric space $M_n$ is universal for a class $\mathcal C$ of proper metric spaces of asymptotic dimension $n$ 
if any metric space in $\mathcal C$ admits a coarse embedding in $M_n$. This is in analogy to the topological dimension theory
where compact spaces of dimension $n$ embed in the Menger compactum $\mu ^n$. It is shown in \cite {DS} that there is no such space
for $n=1$ if we take $\mathcal C$ to be the class of all proper metric spaces of asymptotic dimension $n$. 

On the other hand for hyperbolic spaces Buyalo, Dranishnikov, Schroeder \cite{BDS}
show that if $X$ is visual hyperbolic metric space such that its boundary $\partial X$ is a doubling metric space and $\asdim X=n$ then $X$ quasi-isometrically
embeds in a product of $n+1$-binary metric trees.

As the asymptotic dimension is a `coarse topology' notion it makes sense to consider coarse embeddings of such spaces
rather than quasi-isometries.

%I owe the following example to P. Tselekidis:
%\begin{Ex} 
%Start with the linear tree $L$ with vertices labeled by $\mathbb N$ and glue a copy
%of $L$ to each vertex $2^n$ to obtain a tree $T$. It is easy to see that the growth function of $T$ is linear but there is
%no coarse embedding of $T$ to $\mathbb R$. Similarly one can easily construct examples of graphs with polynomial growth
%$\sim n^k$ which do not coarsely embed in $\mathbb R ^k$.
%\end{Ex}

It is not clear how much bigger dimension one needs in order to achieve a coarse embedding of a graph of growth $\sim r^k$ in $\mathbb R^n$ instead of an 
embedding considered by  Linial, London, Rabinovich \cite{LLR} and  Krauthgamer-Lee. Note that the embeddings of \cite{KL} are weaker than coarse embeddings, as one only requires
that distinct vertices map at distance $\geq 1$ from each other. For example there is an onto embedding in their sense from the linear graph with vertices $\mathbb N$ to the standard Cayley graph of $\mathbb Z^2$-so these embeddings may raise
$\asdim $.

\begin{Qu} Let $G=(V,E)$ be a graph such that its growth function satisfies $\gamma (n)\leq Cn^k$
for some $C>0$. Is there a coarse embedding $f:G\to \mathbb R^{O(k\log k)}$? 
\end{Qu}

%Apparently this is not known even for nilpotent groups of polynomial growth.

One can ask also whether there is a `universal' space for spaces with polynomial growth:

\begin{Qu} Is there a proper metric space $P_k$ such that if $X$ is a metric space of polynomial growth $<Cn^k$
then $X$ coarsely embeds in $P_k$? If so, can one take $P_k$ so that $\asdim P_k=k$ and $P_k$ is of polynomial growth $\sim n^k$?
(it might be necessary here to fix $C$).
\end{Qu}

Manifolds with nonnegative Ricci curvature have been extensively studied (see e.g. \cite{AG}, \cite{Sh}, \cite {Li}).

Our result implies that if $M^n$ is a complete Riemannian manifold of nonnegative Ricci curvature and bounded geometry
then $\asdim M^n\leq n$. One wonders whether the bounded geometry assumption in necessary:

\begin{Qu} Let  $M^n$ be a complete Riemannian manifold of nonnegative Ricci curvature. Is it true that $\asdim M^n\leq n$?
\end{Qu}

%ASK ABOUT COARSE FOLIATIONS?

\end{document}


\begin{thebibliography}{99}

\bibitem{AG} Abresch, U. and Gromoll, D., 1990.{\em On complete manifolds with nonnegative Ricci curvature.} Journal of the American Mathematical Society, 3(2), pp.355-374.

\bibitem{BD}Bell, G. and Dranishnikov, A., 2008. {\em Asymptotic dimension}. Topology and its Applications, 155(12), pp.1265-1296.

\bibitem{BG} I. Benjamini, A. Georgakopoulos 2021, {\em Triangulations of uniform subquadratic growth are quasi-trees}, preprint, arXiv:2106.06443.


\bibitem{BBEGLPS} Bonamy M, Bousquet N, Esperet L, Groenland C, Liu CH, Pirot F, Scott A. {\em Asymptotic Dimension of Minor-Closed Families and Assouad-Nagata Dimension of Surfaces}. arXiv preprint arXiv:2012.02435. 2020 Dec 4, to appear in EJM.

\bibitem{BNLM} Brodskiy, N., Dydak, J., Levin, M. and Mitra, A., 2008.{\em A Hurewicz theorem for the Assouad–Nagata dimension. }Journal of the London Mathematical Society, 77(3), pp.741-756.

\bibitem{BDS} Buyalo, S., Dranishnikov, A. and Schroeder, V., 2007.{\em   Embedding of hyperbolic groups into products of binary trees.} Inventiones mathematicae, 169(1), pp.153-192.


\bibitem{BL}Buyalo, S. and Lebedeva, N., 2008. {\em Dimensions of locally and asymptotically self-similar spaces.} St. Petersburg Mathematical Journal, 19(1), pp.45-65.

\bibitem{DS} Dranishnikov, A.N. and Smith, J., 2007.{\em On asymptotic Assouad–Nagata dimension.} Topology and its Applications, 154(4), pp.934-952.

\bibitem{DZ} Dranishnikov, A. and Zarichnyi, M., 2004.{\em Universal spaces for asymptotic dimension.} Topology and its Applications, 140(2-3), pp.203-225.

\bibitem{DR} Le Donne, E. and Rajala, T., 2015. {\em Assouad dimension, Nagata dimension, and uniformly close metric tangents.} Indiana University Mathematics Journal, pp.21-54.

\bibitem{FP} Fujiwara, K. and Papasoglu, P., 2020. {\em Asymptotic dimension of planes and planar graphs.} arXiv preprint arXiv:2002.01630, to appear in Transactions of AMS.
%\bibitem{FP}
%K. Fujiwara, P. Papasoglu. A quasi-isometric characterization of cacti.
%In preparation. 

%\bibitem{GHT} J.R.Gilbert, J.P. Hutchinson, R.E. Tarjan  {\em A
%separator theorem for graphs of bounded genus } J. Algorithms 5
%(1984) 391-407.

\bibitem{Gr} M.Gromov, {\em Asymptotic invariants of infinite groups}, in Geometric Group Theory, v.
2, 1-295, London Math. Soc. Lecture Note Ser., vol. 182, Cambridge Univ. Press, Cambridge,
1993.

\bibitem{Gr2} Gromov, M., 2007.{\em Metric structures for Riemannian and non-Riemannian spaces.} Springer Science \& Business Media.

\bibitem{He} Heinonen, J., 2012.{\em Lectures on analysis on metric spaces.} Springer Science \& Business Media.




\bibitem{JL}
 J{\o}rgensen M., Lang U.
{\em Geodesic spaces of low Nagata dimension.}
preprint, 2020, 
arXiv:2004.10576, to appear in Ann. Acad. Sci. Fenn. Math.


%\bibitem{KPR} P.Klein,  S.A.Plotkin,  and S. Rao,   {\em Excluded minors, network decomposition, and multicommodity flow}. In Proceedings of the twenty-fifth annual ACM symposium on Theory of computing ( 682-690).1993, June.

\bibitem{KL} Krauthgamer, R. and Lee, J.R., 2007.{\em  The intrinsic dimensionality of graphs.} Combinatorica, 27(5), pp.551-585.



\bibitem{LS} Lang, U. and Schlichenmaier, T., 2005.{\em Nagata dimension, quasisymmetric embeddings, and Lipschitz extensions.} International Mathematics Research Notices, 2005(58), pp.3625-3655.

\bibitem{Lin}Linial, N. {\em Variation on a theme of Levin.} In Open Problems, Workshop on Discrete Metric Spaces and their Algorithmic Applications. 2002.

\bibitem{LLR} Linial, N., London, E. and Rabinovich, Y., 1995. {\em The geometry of graphs and some of its algorithmic applications.} Combinatorica, 15(2), pp.215-245.

\bibitem{Li} Liu, G., 2013. {\em 3-manifolds with nonnegative Ricci curvature.} Inventiones mathematicae, 193(2), pp.367-375.

\bibitem{Y} Yu, G., 1998.{\em The Novikov conjecture for groups with finite asymptotic dimension.} Annals of Mathematics, 147(2), pp.325-355.

%\bibitem{LT} RJ. Lipton , RE. Tarjan .{\em A separator theorem for planar graphs. SIAM Journal on Applied Mathematics}. 1979;36(2):177-89.

%\bibitem{Li} C.H.Liu,  {\em  Asymptotic dimension of minor-closed families and beyond}. preprint 2020.arXiv:2007.08771.


\bibitem{OR}
Mikhail I. Ostrovskii, David Rosenthal, 
{\em Metric dimensions of minor excluded graphs and minor exclusion in groups.} 
Internat. J. Algebra Comput. 25 (2015), no. 4, 541–-554. 

\bibitem{Sh} Shen, Z., 1996. {\em  Complete manifolds with nonnegative Ricci curvature and large volume growth.} Inventiones mathematicae, 125(3), pp.393-404.

\bibitem{ST} \v{S}pakula, J. and Tikuisis, A., 2019.{\em  Relative commutant pictures of Roe algebras.} Communications in Mathematical Physics, 365(3), pp.1019-1048.

\bibitem{T} Tessera, R., 2006.{\em  Asymptotic isoperimetry of balls in metric measure spaces}. Publicacions Matem\'atiques, pp.315-348.

%\bibitem{Ma}
%Manning JF. {\em Geometry of pseudocharacters}, Geometry \& Topology (2005) Jun 14;9(2):1147-85.
\end{thebibliography}
\end{document}
\bye